\documentclass[11pt,A4]{amsart}
\usepackage{ifthen,latexsym,amssymb,amsmath}
\newtheorem{theorem}{Theorem}
\newtheorem{lemma}[theorem]{Lemma}

\begin{document}

\author{Bart\l{}omiej Bzd\c{e}ga}

\address{Adam Mickiewicz University, Pozna\'n, Poland}

\email{exul@amu.edu.pl}



\title{On a certain family of inverse ternary cyclotomic polynomials}

\maketitle

\begin{abstract}
We study a family of inverse ternary cyclotomic polynomials $\Psi_{pqr}$ in which $r\le\varphi(pq)$ is a positive linear combination of $p$ and $q$. We derive a formula for the height of such polynomial and characterize all flat polynomials in this family.
\end{abstract}

\section{Introduction}

Let
$$\Phi_n(x) = \prod_{1\le k \le n,\; (k,n)=1}(x-e^{2k\pi i/n}) = \sum_m a_n(m)x^m$$
be the $n$th cyclotomic polynomial. The $n$th inverse cyclotomic polynomial is defined by the formula
$$\Psi_n(x) = \frac{x^n-1}{\Phi_n(x)} = \sum_m c_n(m)x^m.$$
Like for cyclotomic polynomials, for odd primes $p<q<r$, we say that $\Psi_{pq}$ is binary, $\Psi_{pqr}$ is ternary, etc.

Recall that the height of a given polynomial $F$ is the maximal absolute value of its coefficients. We say that polynomial is flat, if its height equals $1$. Traditionally we denote the height of $\Phi_n$ by $A(n)$ and the height of $\Psi_n$ by $C(n)$.

Ternary inverse cyclotomic polynomials were studied by P.~Moree \cite{Moree-Inverse}. He proved that $C(pqr) \le p-1$ and for every prime $p\ge3$ there are infinitely many pairs $(q,r)$ of primes for which $C(pqr)=p-1$. Additionally he came up with the following bound (\cite{Moree-Inverse}, Theorem 7):
$$C(pqr) \le \max\{\min\{p',q'\},\min\{q-p',p-q'\}\} \text{ for } \deg\Psi_{pqr}<2qr.$$
He also found some flat inverse ternary cyclotomic polynomials.

Let us remark that the case $r>\varphi(pq)=\deg\Phi_{pq}$ is trivial, because by the identity $\Psi_{pqr}(x)=\Psi_{pq}(x^r)\Phi_{pq}(x)$ we have $c_{pqr}(ar+b)=a_{pq}(b)c_{pq}(a)$ for $a\ge0$ and $0\le b < r$. The coefficients of polynomials $\Phi_{pq}$ and $\Psi_{pq}$ are well known, so we can evaluate $c_{pqr}(ar+b)$ easily.

Although there is a substantial research on flat ternary cyclotomic polynomials \cite{Bachman-TernaryFlat, Elder-Flat, Kaplan-TernaryFlat}, we do not know much about flat ternary inverse cyclotomic polynomials. Particularly, no infinite family of such polynomials in which $r\le\varphi(pq)$ was known so far.

In this paper we investigate polynomials $\Psi_{pqr}$ in which $r\le\varphi(pq)$ is a positive linear combination of $p$ and $q$. For this specific type of polynomials we improve some of the results of P.~Moree mentioned above. Our main result is the following theorem.

\begin{theorem} \label{T-Cpqr}
Let $r=\alpha p + \beta q\le\varphi(pq)$, where $\alpha,\beta>0$. Let also $p'\in\{1,2,\ldots,q-1\}$ be the inverse of $p$ modulo $q$ and $q'\in\{1,2,\ldots,p-1\}$ be the inverse of $q$ modulo $p$. Then
$$C(pqr)=\max\left\{\min\left\{\left\lceil\frac{p'}{\alpha}\right\rceil,\left\lceil\frac{q'}{\beta}\right\rceil\right\},\min\left\{\left\lceil\frac{q-p'}{\alpha}\right\rceil,\left\lceil\frac{p-q'}{\beta}\right\rceil\right\}\right\}.$$
\end{theorem}

The above formula is similar to the already mentioned one obtained by P.~Moree. However, our theorem does not require the assumption $\deg\Psi_{pqr}<2qr$. We use Theorem \ref{T-Cpqr} to characterize all flat inverse ternary cyclotomic polynomials $\Psi_{pqr}$ in which $r$ is a positive linear combination of $p$ and $q$.

\begin{theorem}\label{T-flat}
Let $r=\alpha p + \beta q \le \varphi(pq)$, where $\alpha,\beta>0$. Then $\Psi_{pqr}$ is flat if and only if at least one of the following conditions holds:
\begin{itemize}
\item[(a)] $\alpha \ge \max\{p',q-p',\}$,
\item[(b)] $\beta \ge \max\{q',p-q',\}$,
\item[(c)] $\alpha\ge p'$ and $\beta\ge p-q'$,
\item[(d)] $\alpha\ge q-p'$ and $\beta\ge q'$,
\end{itemize}
where $p'$ and $q'$ are like in Theorem \ref{T-Cpqr}.
\end{theorem}

At first one may expect that the set of primes $r\le\varphi(pq)$ satisfying at least one of the conditions (a) -- (d) of Theorem \ref{T-flat} is rather small and consists of primes which are relatively close to $\varphi(pq)$ (if (a) or (b) holds, then $r>pq/2$). Fortunately, the following theorem says that this is not the truth in general.

\begin{theorem} \label{T-flatcond}
Let $S(p,q)$ denote the set of primes $r\le\varphi(pq)$ of form $\alpha p + \beta q$, $\alpha,\beta >0$, for which $\Psi_{pqr}$ is flat. Then
\begin{itemize}
\item[(i)] For every $M>0$ and every $p$ there exists such $q$ that $\# S(p,q)>M$.
\item[(ii)] For every $\varepsilon>0$ there exists a triple $(p,q,r)$ of primes such that $r\in S(p,q)$ and $r<\varepsilon\varphi(pq)$.
\end{itemize}
\end{theorem}
So we reveal a new, vast family of nontrivially flat ternary inverse cyclotomic polynomials.

Our paper is organized in the following way. In section \ref{S-coeff} we derive a formula for $c_{pqr}(m)$, the $m$th coefficient of $\Psi_{pqr}$. In section \ref{S-height} we prove Theorem \ref{T-Cpqr}. Finally, in section \ref{S-flat} we prove Theorems \ref{T-flat} and \ref{T-flatcond}.

\section{Preliminaries}

In this section we recall some basic properties of standard and inverse cyclotomic polynomials. All of them can be found in \cite{LamLeung-Binary} or \cite{Moree-Inverse}. Let us start with the properties of their degrees:
$$\deg \Phi_n = \varphi(n), \qquad \deg\Psi_n = n-\varphi(n).$$
Particularly $\deg\Psi_{pqr}=qr+rp+pq-p-q-r+1$.

We say that a polynomial $F$ is reciprocal if $F(x)=x^{\deg F}F(1/x)$ and anti-reciprocal if $F(x)=-x^{\deg F}F(1/x)$. It is known that all cyclotomic polynomials except of $\Phi_1$ are reciprocal and all inverse cyclotomic polynomials except of $\Psi_1$ are anti-reciprocal.

In our investigations we need to know the coefficients of $\Phi_{pq}$. The following lemma, proved in \cite{LamLeung-Binary}, derives a formula on $a_{pq}(m)$.

\begin{lemma} \label{L-Phipq}
Let $m\in\{0,1,\ldots,pq-1\}$ and let $u,v$ be the unique numbers such that $m\equiv up+vq \pmod{pq}$ and $0 \le u < q$, $0 \le v < p$. Then we have
$$a_{pq}(m) = \left\{\begin{array}{ll}
1, & \text{if } u < p' \text{ and } v < q', \\
-1, & \text{if } u \ge p' \text{ and } v \ge q', \\
0, & \text{otherwise}.
\end{array}\right.$$
\end{lemma}

Following P.~Moree \cite{Moree-Inverse}, we define the polynomial
$$f_{pqr}(x) = (1+x^r+\ldots+x^{(p-1)r})\Phi_{pq}(x) = \sum_m e_{pqr}(m)x^m.$$
Let also
$$\tau(pqr) = \deg f_{pqr} = (p-1)r+\varphi(pq)=(p-1)(q+r-1).$$
The next lemma was proved in \cite{Moree-Inverse}.

\begin{lemma} \label{L-Psipqr}
For primes $p<q<r$ we have
$$\Psi_{pqr}(x)=(x^{qr}-1)f_{pqr}(x).$$
\end{lemma}

By this lemma $c_{pqr}(m)=e_{pqr}(m-qr)-e_{pqr}(m)$, so if we want to determine the coefficients of $\Psi_{pqr}$, we need to know the coefficients of $f_{pqr}$.

Let us remark that in the formula
$$\Psi_{pqr}(x) = -\frac{(1-x)(1-x^{qr})(1-x^{rp})(1-x^{pq})}{(1-x^p)(1-x^q)(1-x^r)}$$
one can replace the assumption that $p$, $q$, $r$ are primes by the assumption that they are pairwise coprime. This way we receive the definition of inverse inclusion-exclusion polynomial $\Psi_{p,q,r}$. Theorems \ref{T-Cpqr} and \ref{T-flat} hold also for inverse inclusion-exclusion polynomials and they can be proved by analogous methods.

\section{Coefficients of $f_{pqr}$} \label{S-coeff}

The following lemma, partially proved in \cite{Moree-Inverse}, derives a formula on coefficients of $f_{pqr}$ in terms of coefficients of $\Phi_{pq}$. We remark that this is true for all primes $r$.

\begin{lemma} \label{L-fcoeff}
The following equalities hold:
\begin{itemize}
\item[(i)] if $m<pr$, then $e_{pqr}(m)=\sum_{j=0}^{\lfloor m/r \rfloor}a_{pq}(m-jr)$,
\item[(ii)] for $pq<m<pr$ we have $e_{pqr}(m)=e_{pqr}(m-r)$.
\item[(iii)] if $pr \le m \le \tau(pqr)$, then $e_{pqr}(m)=e_{pqr}(m')$, where $m'=\tau(pqr)-m<pr$.
\end{itemize}
\end{lemma}

\begin{proof}
Case (i) follows directly from the definition of $f_{pqr}$. The polynomial $f_{pqr}$ is reciprocal as a product of reciprocal polynomials, so $e_{pqr}(m)=e_{pqr}(\tau(pqr)-m)$. Because $\tau(pqr)<2pr$, for $m\ge pr$ we have $m'=\tau(pqr)-m<pr$, so (iii) holds. To prove (ii) we observe that for $pq<m<pr$ we have $a_{pq}(m)=0$ and then by (i)
$$e_{pqr}(m) = a_{pq}(m) + \sum_{j=1}^{\lfloor m/r \rfloor}a_{pq}(m-jr) = \sum_{j=0}^{\lfloor (m-r)/r \rfloor}a_{pq}(m-r-jr) = e_{pqr}(m-r),$$
which completes the proof.
\end{proof}

Now we use Lemmas \ref{L-Phipq} and \ref{L-fcoeff} to determine coefficients of $f_{pqr}$. We do it for the exponents not greater than $pq$, since for greater ones we can use (ii) and (iii) of Lemma \ref{L-fcoeff}. In order to simplify the notation, for a finite set $A$ we define
$$\min_{\ge0} A = \max\{0,\min A\}.$$

\begin{theorem} \label{T-fcoeff}
Let $r=\alpha p + \beta q \le \varphi(pq)$, $\alpha,\beta>0$ and $m<pq$. Put $m=(a-1)r+b$, where $0\le b < r$.
\begin{itemize}
\item[(i)] If $b=up+vq$ for some integers $0\le u < q$ and $0\le v < p$, then
$$e_{pqr}(m)=\min_{\ge0}\left\{a,\left\lceil\frac{p'-u}{\alpha}\right\rceil,\left\lceil\frac{q'-v}{\beta}\right\rceil\right\}.$$
\item[(ii)] If $b+pq=up+vq$ for some integers $0\le u < q$ and $0\le v < p$, then we define
$$j_0=\min\left\{\left\lceil\frac{q-u}{\alpha}\right\rceil,\left\lceil\frac{p-v}{\beta}\right\rceil\right\}, \qquad a^*=a-j_0,$$
$$(u^*,v^*) = \left\{\begin{array}{ll}
(u+j_0\alpha-q,v+j_0\beta) & \text {if } j_0=\lceil(q-u)/\alpha\rceil, \\
(u+j_0\alpha,v+j_0\beta-p) & \text {if } j_0=\lceil(p-v)/\beta\rceil.
\end{array}\right.$$
We have $e_{pqr}(m) = e^+_{pqr}(m) - e^-_{pqr}(m)$, where
\begin{align*}
e^+_{pqr}(m) & = \min_{\ge0}\left\{a^*,\left\lceil\frac{p'-u^*}{\alpha}\right\rceil,\left\lceil\frac{q'-v^*}{\beta}\right\rceil\right\}, \\
e^-_{pqr}(m) & = \min_{\ge0}\left\{\min\left\{a,\left\lceil\frac{q-u}{\alpha}\right\rceil,\left\lceil\frac{p-v}{\beta}\right\rceil\right\}\right. \\
& \quad -\left.\max\left\{0,\left\lceil\frac{p'-u}{\alpha}\right\rceil,\left\lceil\frac{q'-v}{\beta}\right\rceil\right\}\right\}.
\end{align*}
\end{itemize}
\end{theorem}

One can easily prove that every $0\le b < pq$ can be written in exactly one of forms: $up+vq$ or $up+vq-pq$, where $0\le u<q$ and $0\le v<p$. So cases (i) and (ii) cover all possible values of $b$.

\begin{proof}
First consider $b=up+vq$. Then for $0\le j < a$ we have
$$b+jr = (u+j\alpha)p + (v+j\beta)q.$$
Notice that $b+jr \le b+(a-1)r=m<pq$, so $u+j\alpha<q$ and $v+j\beta<p$. By Lemma \ref{L-Phipq}
$$a_{pq}(b+jr) = \left\{\begin{array}{ll}
1, & \text{if } u+j\alpha < p' \text{ and } v+j\beta < q', \\
0, & \text{otherwise},
\end{array}\right.$$
where the case $a_{pq}(b+jr)=-1$ was omitted, as the inequalities $u+j\alpha \ge p'$ and $v+j\beta \ge q'$ cannot hold at the same time (if they both held then we would have $pq > b+jr \ge pp'+qq'=pq+1$, a contradiction). So for $0 \le j < a$ we have $a_{pq}(b+jr)=1$ if and only if both equalities
$$0\le j < \frac{p'-u}{\alpha}, \qquad 0\le j < \frac{q'-v}{\beta}$$
hold. Thus by Lemma \ref{L-fcoeff}
$$e_{pqr}(m) = \sum_{j=0}^{a-1}a_{pq}(b+jr) = \min_{\ge0}\left\{a,\left\lceil\frac{p'-u}{\alpha}\right\rceil,\left\lceil\frac{q'-v}{\beta}\right\rceil\right\}$$
as desired. We have proved (i).

Now consider $b+pq=up+vq$. For $0\le j < a$ we have
$$b+jr = (u+j\alpha)p + (v+j\beta)q - pq.$$
As long as $u+j\alpha<q$ and $v+j\beta<p$ (or equivalently $j<j_0$), $b+jr$ is clearly not a nonnegative linear combination of $p$ and $q$, while for $j\ge j_0$ it is. Therefore we put
$$-h^- = \sum_{j=0}^{\min\{j_0,a\}-1}a_{pq}(b+jr), \qquad h^+ = \sum_{j=\min\{j_0,a\}}^{a-1}a_{pq}(b+jr).$$
(if a sum is empty then it equals $0$). By Lemma \ref{L-fcoeff} we have $e_{pqr}(m) = h^+-h^-$, so we need to prove that $h^-=e^-_{pqr}(m)$ and $h^+=e^+_{pqr}(m)$.

The value of $h^-$ is easier to determine. For $j<j_0$ the inequalities $u+j\alpha < p'$ and $v+j\beta < q'$ cannot hold at the same time, because if $kp+lq \ge pq$, then $k\ge p'$ or $l\ge q'$. So for $j<j_0$ we receive
$$a_{pq}(b+jr) = \left\{\begin{array}{ll}
-1, & \text{if } u+j\alpha \ge p' \text{ and } v+j\beta \ge q', \\
0, & \text{otherwise}.
\end{array}\right.$$
Therefore for $0\le j < \min\{j_0,a\}$ we have $a_{pq}(b+jr)=-1$ if and only if both inequalities
$$\frac{p'-u}{\alpha} \le j < \frac{q-u}{\alpha}, \qquad \frac{q'-v}{\beta} \le j < \frac{p-v}{\beta}$$
hold. Thus
\begin{align*}
h^- & = \min_{\ge0}\left\{\min\left\{j_0,a,\left\lceil\frac{q-u}{\alpha}\right\rceil,\left\lceil\frac{p-v}{\beta}\right\rceil\right\}\right. \\
& \quad - \left.\max\left\{0,\left\lceil\frac{p'-u}{\alpha}\right\rceil,\left\lceil\frac{q'-v}{\beta}\right\rceil\right\}\right\} \\
& = e^-_{pqr}(m),
\end{align*}
by of the definition of $j_0$.

Now we determine $h^+$. If $j_0\ge a$, then $h^+=0$, so further we assume that $j_0<a$. We have two analogous cases here: $j_0=\lceil(q-u)/\alpha\rceil$ and $j_0=\lceil(p-v)/\beta\rceil$. Let us consider the first one. For $j_0 \le j < a$ we have
$$b+jr = (u+j\alpha-q)p + (v+j\beta)q$$
with $0\le u+j\alpha-q < q$ and $v+j\beta < p$. We have $u^* = u+j_0\alpha-q$ and $v^* = v+j_0\beta$. Additionally put
$$b^* = b+j_0r = (u+j_0\alpha-q)p + (v+j_0\beta)q = u^*p+v^*q.$$
Clearly $b^*\le b+(a-1)r<pq$, so $u^*<q$ and $v^*<p$. Note that $b^*$ does not have to be smaller than $r$, however, we can still use the arguments from the proof of (i) to obtain that
$$h^+ = \sum_{j=0}^{a^*-1}a_{pq}(b^*+jr) = e^+_{pqr}(m).$$
This formula remains correct even if $j_0\ge a$. Applying the analogous argument to the case $j_0=\lceil(p-v)/\beta\rceil$ we complete the proof.
\end{proof}

From the results of this section, we obtained an algorithm which instantly computes the value of $c_{pqr}(k)$. By Lemma \ref{L-Psipqr} we have $c_{pqr}(k) = e_{pqr}(k-qr) - e_{pqr}(k)$. Then using Lemma \ref{L-fcoeff} we reduce computing $c_{pqr}(k)$ to computing $e_{pqr}(m_1)$ and $e_{pqr}(m_2)$ for some $m_1,m_2 < pq$. Finally we apply Theorem \ref{T-fcoeff} to evaluate $e_{pqr}(m_1)$ and $e_{pqr}(m_2)$.

\section{The height of $f_{pqr}$ and $\Psi_{pqr}$} \label{S-height}

In this section we evaluate the height of $f_{pqr}$ and compare it with the height of $\Psi_{pqr}$.

\begin{lemma} \label{L-Hf}
Let $H(pqr)$ denotes the height of $f_{pqr}$. Then for $r=\alpha p + \beta q \le \varphi(pq)$, $\alpha,\beta>0$, we have
$$H(pqr) = \max\left\{\min\left\{\left\lceil\frac{p'}{\alpha}\right\rceil,\left\lceil\frac{q'}{\beta}\right\rceil\right\},\min\left\{\left\lceil\frac{q-p'}{\alpha}\right\rceil,\left\lceil\frac{p-q'}{\beta}\right\rceil\right\}\right\}.$$
\end{lemma}

\begin{proof}
By Lemma \ref{L-fcoeff} we can restrict our considerations to $m<pq$ and use Theorem \ref{T-fcoeff}. We receive the inequalities
\begin{align*}
e_{pqr}(m) & \le \min\left\{\left\lceil\frac{p'}{\alpha}\right\rceil,\left\lceil\frac{q'}{\beta}\right\rceil\right\}, \\
-e_{pqr}(m) & \le \max_{u,v}\left\{\min\left\{\left\lceil\frac{q-u}{\alpha}\right\rceil,\left\lceil\frac{p-v}{\beta}\right\rceil\right\}-\max\left\{\left\lceil\frac{p'-u}{\alpha}\right\rceil,\left\lceil\frac{q'-v}{\beta}\right\rceil\right\}\right\} \\
& \le \min\left\{\max_u\left\{\left\lceil\frac{q-u}{\alpha}\right\rceil-\left\lceil\frac{p'-u}{\alpha}\right\rceil\right\},\max_v\left\{\left\lceil\frac{p-v}{\beta}\right\rceil-\left\lceil\frac{q'-v}{\beta}\right\rceil\right\}\right\} \\
& = \min\left\{\left\lceil\frac{q-p'}{\alpha}\right\rceil,\left\lceil\frac{p-q'}{\beta}\right\rceil\right\}.
\end{align*}

To complete the proof we will show that we have equalities: in the first inequality for $m=m_1$ and in the second one for $m=m_2$, where
\begin{align*}
m_1 & = \left(\min\left\{\left\lceil\frac{p'}{\alpha}\right\rceil,\left\lceil\frac{q'}{\beta}\right\rceil\right\}-1\right)r, \\
m_2 & = \left(\min\left\{\left\lceil\frac{q-p'}{\alpha}\right\rceil,\left\lceil\frac{p-q'}{\beta}\right\rceil\right\}-1\right)r+1.
\end{align*}
In order to use Theorem \ref{T-fcoeff}, we have to show that $m_1,m_2<pq$. Indeed,
\begin{align*}
m_1 & \le \left(\left\lceil\frac{p'}{\alpha}\right\rceil-1\right)\alpha p + \left(\left\lceil\frac{q'}{\beta}\right\rceil-1\right)\beta q \le (pp'-1) + (qq'-1) = pq-1, \\
m_2 & \le \left(\left\lceil\frac{q-p'}{\alpha}\right\rceil-1\right)\alpha p + \left(\left\lceil\frac{p-q'}{\beta}\right\rceil-1\right)\beta q \\
& \le (pq-pp'-1) + (pq-qq'-1) = pq-1.
\end{align*}
For $m=m_1$ we have $u=v=0$, so the first desired equality follows easily from (i) of Theorem \ref{T-fcoeff}. For $m=m_2$ we have $u=p'$, $v=q'$ and $a=j_0$, $e^+_{pqr}(m)=0$, so the second desired equality holds by (ii) of Theorem \ref{T-fcoeff}. Thus the proof is done.
\end{proof}

Now we are ready to prove our main result.

\begin{proof}[Proof of Theorem \ref{T-Cpqr}]
By Lemma \ref{L-Hf}, we have to prove that $C(pqr)=H(pqr)$.

First let us verify that $C(pqr)\ge H(pqr)$. Because $f_{pqr}$ is reciprocal, $|e_{pqr}(m)|=H(pqr)$ for some $m\le\tau(pqr)/2$. Since $\tau(pqr)<2qr$, we have $m<qr$ and hence $e_{pqr}(m-qr)=0$. Thus by Lemma \ref{L-Psipqr}
$$C(pqr) \ge |c_{pqr}(m)| = |e_{pqr}(m-qr)-e_{pqr}(m)| = |-e_{pqr}(m)| = H(pqr).$$

The opposite inequality is much harder to prove. By Lemma \ref{L-Psipqr} we have $|c_{pqr}(m)| = |-e_{pqr}(m)| \le H(pqr)$ for $m<qr$. By the anti-reciprocity of $\Psi_{pqr}$ we also have $|c_{pqr}(m)|\le H(pqr)$ for $m>\tau(pqr)$. Therefore we can restrict our considerations to $qr \le m \le \tau(pqr)$. In this case
\begin{align*}
c_{pqr}(m) & = e_{qpr}(m-qr)-e_{pqr}(m) = e_{pqr}(m-qr) - e_{pqr}(\tau(pqr)-m) \\
& = e_{pqr}(m_1) - e_{pqr}(m_2),
\end{align*}
where $m_1=m-qr$ and $m_2=\tau(pqr)-m$. Additionally,
$$m_1+m_2 = \tau(pqr)-qr = \varphi(pq) - (q-1-p)r < pq.$$
Hence $0\le m_1 < pq$, $0\le m_2 < pq$ and so we can use Theorem \ref{T-fcoeff}. We will show that $e_{pqr}(m_1)$ and $e_{pqr}(m_2)$ cannot have opposite signs, which actually completes the proof. Without loss of generality we can assume that $e_{pqr}(m_1)>0$.

For $i\in\{1,2\}$ put
$$m_i=(a_i-1)r+b_i, \qquad 0\le b_i<r,$$
$$b_i\equiv u_ip+v_iq \pmod{pq}, \qquad 0\le u_i < q, \qquad 0\le v_i < p.$$

We have
$$m_2 = \tau(pqr)-qr-m_1 = (p-q-a_1)r + \varphi(pq)-b_1 = (t+p-q-a_1)r+b_2,$$
where $\varphi(pq)-b_1=tr+b_2$ with $0\le b_2 < r$. Then $a_2 = t+p-q-a_1+1$. Now we consider some cases, in which we determine different values of $u_2$, $v_2$.

\medskip\noindent\texttt{Case (1):} $b_1=u_1p+v_1q$. Here
\begin{align*}
b_2 & = \varphi(pq)-b_1-tr = (p'-1)p+(q'-1)q - (u_1p+v_1q) - t(\alpha p + \beta q) \\
& = (p'-1-u_1-t\alpha)p + (q'-1-v_1-t\beta)q,
\end{align*}
so $u_2\equiv p'-1-u_1-t\alpha \pmod q$ and $v_2\equiv q'-1-v_1-t\beta \pmod p$.

Both numbers $p'-1-u_1-t\alpha$ and $q'-1-v_1-t\beta$ cannot be negative at the same time, since $b_2\ge0$. If both are positive, then they equal $u_2$ and $v_2$ and $e_{pqr}(m_2)\ge0$ by Theorem \ref{T-fcoeff}. Therefore we have to consider the situation in which one of these numbers is negative and one is positive. Without loss of generality, we assume that $p'-1-u_1-t\alpha<0$. Then
$$u_2 = q+p'-1-u_1-t\alpha, \qquad v_2 = q'-1-v_1-t\beta, \qquad b_2 = u_2p+v_2q-pq.$$

\noindent\texttt{Case (2):} $b_1=u_1p+v_1q-pq$. Here $u_1\ge p'$ or $v_1\ge q'$. Without loss of generality we assume that $u_1\ge p'$. Then
\begin{align*}
b_2 & = \varphi(pq)-b_1-tr = (p'-1)p+(q'-1)q - (u_1p+v_1q) + pq - t(\alpha p + \beta q) \\
& = (q+p'-1-u_1-t\alpha)p + (q'-1-v_1-t\beta)q.
\end{align*}
By the similar argument to one used in case (1), we consider two subcases in which the signs of $q+p'-1-u_1-t\alpha$ and $q'-1-v_1-t\beta$ are opposite:

\medskip\noindent\texttt{Case (2a):} $u_2 = 2q+p'-1-u_1-t\alpha$, $v_2 = q'-1-v_1-t\beta$, $b_2 = u_2p+v_2q-pq$,

\medskip\noindent\texttt{Case (2b):} $u_2 = q+p'-1-u_1-t\alpha$, $v_2 = p+q'-1-v_1-t\beta$, $b_2 = u_2p+v_2q-pq$.

\medskip Now we show that in cases (1) and (2a) we have $e_{pqr}(m_2)\ge0$. Note that $b_2 = u_2p+v_2q-pq$ and $v_2=q'-1-v_1-t\beta$ in both these cases, so we estimate
\begin{align*}
\max\left\{0,\left\lceil\frac{p'-u_2}{\alpha}\right\rceil,\left\lceil\frac{q'-v_2}{\beta}\right\rceil\right\} & \ge \left\lceil\frac{q'-v_2}{\beta}\right\rceil = \left\lceil\frac{v_1+1}{\beta}\right\rceil + t \ge t, \\
\min\left\{a_2,\left\lceil\frac{q-u_2}{\alpha}\right\rceil,\left\lceil\frac{p-v_2}{\beta}\right\rceil\right\} & \le a_2 = t+p-q-a_1+1 \le t.
\end{align*}
Hence $e_{pqr}(m_2)\ge e^-_{pqr}(m_2)=0$ by Theorem \ref{T-fcoeff}, which completes the proof in cases (1) and (2a).

\medskip It remains to prove that $e_2(m)\ge0$ in case (2b). We will use variables $u^*_1$, $v^*_1$, $a^*_1$, $j_{0,1}$, which we define like in part (ii) of Theorem \ref{T-fcoeff}. As we assumed $e_{pqr}(m_1)>0$, by Theorem \ref{T-fcoeff} we have $e^+_{pqr}(m_1)>0$. Thus $u^*_1<p'$ and $v^*_1<q'$. So $u^*_1 \neq u_1+j_{0,1}\alpha$ because $u_1\ge p'$. Therefore $u^*_1=u_1+j_{0,1}\alpha-q$ and $v^*_1=v_1+j_{0,1}\beta<q'$, which implies $v_1<q'$ and $j_{0,1} = \lceil(q-u_1)/\alpha\rceil$.

Let us assume that $e_{pqr}(m_2)<0$. We will show that it leads to a contradiction. We have then $e^+_{pqr}(m_1),e^-_{pqr}(m_2)>0$ and so
\begin{align*}
e^+_{pqr}(m_1) & = \min\left\{a^*_1,\left\lceil\frac{p'-u^*_1}{\alpha}\right\rceil,\left\lceil\frac{q'-v^*_1}{\beta}\right\rceil\right\} \le \left\lceil\frac{q'-v^*_1}{\beta}\right\rceil \\
& = \left\lceil\frac{q'-v_1}{\beta}\right\rceil-\left\lceil\frac{q-u_1}{\alpha}\right\rceil.
\end{align*}
Additionally,
\begin{align*}
e^-_{pqr}(m_2) & = \min\left\{a_2,\left\lceil\frac{q-u_2}{\alpha}\right\rceil,\left\lceil\frac{p-v_2}{\beta}\right\rceil\right\} \\
& \quad - \max\left\{0,\left\lceil\frac{p'-u_2}{\alpha}\right\rceil,\left\lceil\frac{q'-v_2}{\beta}\right\rceil\right\} \\
& \le \left\lceil\frac{p-v_2}{\beta}\right\rceil - \left\lceil\frac{p'-u_2}{\alpha}\right\rceil = \left\lceil\frac{v_1+1-q'}{\beta}\right\rceil - \left\lceil\frac{u_1+1-q}{\alpha}\right\rceil.
\end{align*}
Let $l_1 = q'-v_1$ and $l_2=q-u_1$. Combining the above bounds we conclude that
$$e^+_{pqr}(m_1) + e^-_{pqr}(m_2) \le \left(\left\lceil\frac{l_1}{\beta}\right\rceil+\left\lceil\frac{1-l_1}{\beta}\right\rceil\right) - \left(\left\lceil\frac{l_2}{\alpha}\right\rceil+\left\lceil\frac{1-l_2}{\alpha}\right\rceil\right).$$
Since for all positive integers  $l$ and $\gamma$ we have
$$\left(\left\lceil\frac{l}{\gamma}\right\rceil+\left\lceil\frac{1-l}{\gamma}\right\rceil\right)\in\{1,2\},$$
we receive $e^+_{pqr}(m_1) + e^-_{pqr}(m_2)\le1$, contradicting $e^+_{pqr}(m_1),e^-_{pqr}(m_2)>0$. The proof is finally completed.
\end{proof}

\section{Flat polynomials} \label{S-flat}

\begin{proof}[Proof of Theorem \ref{T-flat}]
By Theorem \ref{T-Cpqr}, the polynomial $\Psi_{pqr}$ is flat if and only if
$$\max\left\{\min\left\{\left\lceil\frac{p'}{\alpha}\right\rceil,\left\lceil\frac{q'}{\beta}\right\rceil\right\},\min\left\{\left\lceil\frac{q-p'}{\alpha}\right\rceil,\left\lceil\frac{p-q'}{\beta}\right\rceil\right\}\right\}=1.$$
This is equivalent to
$$(\alpha\ge p' \text{ or } \beta \ge q') \text{ and } (\alpha \ge q-p' \text{ or } \beta \ge p-q'),$$
which by the logical distributive laws is equivalent to
\begin{center}
(a) or (b) or (c) or (d)
\end{center}
as desired.
\end{proof}

\begin{proof}[Proof of Theorem \ref{T-flatcond}]
Put $q=tp+1$. Then $p'=q-t$ and $q'=1$, so (d) transforms into the condition
$$\alpha \ge t \text{ and } \beta \ge 1.$$
Therefore all primes $r$ from the arithmetic progression
$$(q+tp, q+(t+1)p, \ldots) = (2tp+1, (2t+1)p+1, \ldots)$$
satisfy the condition (d). Recall that $\pi(x;a,n)$ denotes the number of primes $r\le x$ satisfying $r\equiv a \pmod n$. By the Dirichlet's theorem on primes in arithmetic progressions and by the above observations, we have
\begin{align*}
\# S(p,q) & \ge \pi(\varphi(pq);1,p) - \pi(2tp;1,p) \\
& = \pi((p-1)tp;1,p) - \pi(2tp;1,p) \\
& \sim \frac1{p-1}\left(\frac{(p-1)tp}{\log t + \log(p^2-p)}-\frac{2tp}{\log t + \log(2p)}\right) \\
& \sim \frac{p(p-3)}{p-1}\frac{t}{\log t} \to \infty
\end{align*}
with $t\to\infty$, as we assumed $p>3$ Hence (1) is proved.

To prove (ii) we again put $q=tp+1$ and we assume that $p>5$, so $4tp+1=4(q-1)+1<\varphi(pq)$. Then once more we use the Dirichlet's theorem to show that the arithmetic progression
$$(2tp+1,(2t+1)p+1,\ldots,4tp+1)$$
contains asymptotically $\frac{2p}{p-1}\frac{t}{\log t}$ primes as $t\to\infty$, and hence for $t$ large enough it contains at least one prime. If $r$ is a prime contained in this progression, then $r\in S(p,q)$ and
$$\frac{r}{\varphi(pq)} \le \frac{4tp+1}{(p-1)tp} < \frac5p$$
Therefore for every prime $p>5$ there exist a prime $q$ and a prime $r\in S(p,q)$ such that $r/\varphi(pq)<5/p$. We can chose $p$ arbitrarily, so the proof is completed.
\end{proof}

\section*{Acknowledgments}

The author is partially supported by NCN grant. He would like to thank Bartosz Naskr\c{e}cki for his impressive and inspiring numerical simulations.

\end{document}